\documentclass[a4paper,10pt]{article}
\usepackage[utf8x]{inputenc}
\usepackage{epsfig}
\usepackage{amssymb}
\usepackage{amsmath}
\usepackage{amsfonts}
\usepackage{moreverb}
\usepackage{graphicx}
\usepackage{rotating}
\usepackage{multirow}

\usepackage[colorlinks,bookmarksopen,bookmarksnumbered,citecolor=red,urlcolor=red]{hyperref}
\title{A well-balanced finite volume scheme for 1D hemodynamic simulations\footnote{We thank ENDOCOM ANR for financial
 support, we thank Fran\c{c}ois Bouchut for fruitful discussions.}}
\author{{Olivier Delestre$^{\dag,}$}\footnote{CNRS \& UPMC Universit\'e Paris 06, UMR 7190, 4 place Jussieu, Institut Jean
 le Rond d'Alembert, Bo\^ite 162, F-75005 Paris, France ; delestre@ida.upmc.fr}\footnote{Presently at:
 Laboratoire J.A. Dieudonn\'e \& EPU Nice-Sophia, Universit\'e de Nice Sophia Antipolis, France ;
 delestre@polytech.unice.fr}
 and {Pierre-Yves Lagr\'ee}\footnote{CNRS \& UPMC Universit\'e Paris 06, UMR 7190, 4 place Jussieu, Institut Jean le
 Rond d'Alembert, Bo\^ite 162, F-75005 Paris, France ; pierre-yves.lagree@upmc.fr}}

\begin{document}

\maketitle

\begin{abstract} English version:
 We are interested in simulating blood flow in arteries with variable elasticity with a one dimensional
 model. We present a well-balanced finite volume scheme based on the recent developments in shallow water equations
 context.
 We thus get a mass conservative scheme which also preserves equilibria of $Q=0$. This numerical method is tested on
 analytical tests.

Version Fran\c{c}aise : Nous nous int\'eressons \`a la simulation d'\'ecoulements sanguins dans des art\`eres dont les parois
 sont \`a \'elasticit\'e variable. Ceci est mod\'elis\'e \`a l'aide d'un mod\`ele unidimensionnel. Nous pr\'esentons
 un sch\'ema "volume fini \'equilibr\'e" bas\'e sur les d\'eveloppements r\'ecents effectu\'es pour la r\'esolution du
 syst\`eme de Saint-Venant. Ainsi, nous obtenons un sch\'ema qui pr\'eserve le volume de fluide ainsi que les
 \'equilibres au repos: $Q=0$. Le sch\'ema introduit est test\'e sur des solutions analytiques.
 \end{abstract}
%
%
%
\section*{Introduction}
We consider the following system of mass and momentum conservation with non dimensionless parameters and variables,
 which is the 1D model of blood flow in an artery or a vessel with non uniform elasticity (it is rewritten in a
 conservative form compared to what we usually find in litterature)
\begin{equation}
 \left\{\begin{array}{l}
         \partial_t A+\partial_x Q=0\\
	\partial_t Q+\partial_x \left[\dfrac{Q^2}{A}+\dfrac{1}{3\sqrt{\pi}\rho}kA^{3/2}\right]
=\dfrac{A}{\sqrt{\pi}\rho}\left(\partial_x {\cal A}_0-\dfrac{2}{3}\sqrt{A}\partial_x k\right)-C_f\dfrac{Q}{A}
        \end{array} \right.,\label{eq:system1}
\end{equation}
with ${\cal A}_0=k\sqrt{A_0}$ and where $A(x,t)$ is the cross-section area ($A=\pi R^2$ with $R$ the radius of the
 arteria), $Q(x,t)=A(x,t)u(x,t)$ the flow rate or the discharge, $u(t,x)$ the mean flow velocity, $\rho$ the blood
 density, $A_0(x)$ the cross section at rest and $k(x)$ the stiffness of the artery. System \eqref{eq:system1} is into
 the form of the Saint-Venant problem with variable pressure presented in \cite{Bouchut04}. We have to mention that
 arterial pulse wavelengths are long enough to justify the use of a 1D model rather than a 3D model when a global
 simulation of blood flow in the cardiovascular system is needed.

\section{Numerical method}
 Since \cite{Bermudez94,Greenberg96}, it is well known (in the shallow water community) that the scheme should be
 well-balanced for good source term treatment, {\it i.e.} the scheme should preserve at least some steady states.
 For system \eqref{eq:system1}, we should preserve at least the "man at eternal rest" or "dead man equilibrium"
 \cite{Delestre11} (without artifacts such as \cite{Kirkman03}), it writes
\begin{equation}
\left\{\begin{array}{l}
  u=0\\
  \dfrac{1}{\sqrt{\pi} \rho}A^{3/2}k-\dfrac{A}{\sqrt{\pi} \rho}{\cal A}_0=Cst        
       \end{array}\right.,
\end{equation}
this means that steady states at rest are preserved (this is the analogous of the "lake at rest" equilibrium). Thus we
 use the scheme proposed in \cite[p.93-94]{Bouchut04} for that kind of model. This is a finite volume scheme with a
 modification of the hydrostatic reconstruction (introduced in \cite{Audusse04c,Bouchut04} for the shallow water
 model).

\subsection{Convective step}
For the homogeneous system
\begin{equation}
\partial_t U+\partial_x F(U,Z)=0, 
\end{equation}
which is \eqref{eq:system1} with:
\begin{equation*}
 U=\left(\begin{array}{c}
          A\\
	  Q
         \end{array}\right),\;
Z=\left(\begin{array}{c}
         {\cal A}_0\\
	k
        \end{array}\right)
\;\text{and}\;
F(U,k)=\left(\begin{array}{c}
              Q\\
	    Q^2/A+kA^{3/2}/(3\sqrt{\pi}\rho)
             \end{array}\right),
\end{equation*}
an explicit first order conservative scheme writes
\begin{equation}
 \dfrac{U_i^{n+1}-U_i^n}{\Delta t}+\dfrac{F_{i+1/2}^n-F_{i-1/2}^n}{\Delta x}=0,\label{eq:scheme-hom}
\end{equation}
where $U_i^n$ is an approximation of $U$
\begin{equation*}
 U_i^n\simeq \dfrac{1}{\Delta x} \int_{x_{i-1/2}}^{x_{i+1/2}} U(x,t_n)dx.
\end{equation*}
$i$ refers to the cell $C_i=(x_{i-1/2},x_{i+1/2})=(x_{i-1/2},x_{i-1/2}+\Delta x)$ and $n$ to time $t_n$ with
 $t_{n+1}-t_n=\Delta t$.\\
The two points numerical flux
\begin{equation*}
F_{i+1/2}^n={\cal F}(U_i^n,U_{i+1}^n,k_{i+1/2}^*), 
\end{equation*}
 with $k_{i+1/2}^*=\max(k_i,k_{i+1})$, is an approximation of the flux function $F(U,Z)$ at the cell interface $i+1/2$. This
 numerical flux will be detailled in subsection \ref{HLL}.

\subsection{Source terms treatment}
In system \eqref{eq:system1}, the terms $A(\partial_x {\cal A}_0-2\sqrt{A}\partial_x k/3)/(\sqrt{\pi}\rho)$ are
 involved in steady states preservation, they need a well-balanced treatment: the variables are reconstructed locally
 thanks to a variant of the hydrostatic reconstruction \cite[p.93-94]{Bouchut04}
\begin{equation}
 \left\{\begin{array}{l}
         \sqrt{A_{i+1/2L}}=\max(k_i\sqrt{A_i}+\min(\Delta {{\cal A}_0}_{i+1/2},0),0)/k^*_{i+1/2}\\
	  U_{i+1/2L}=(A_{i+1/2L},A_{i+1/2L}.u_i)^t\\
	  \sqrt{A_{i+1/2R}}=\max(k_{i+1}\sqrt{A_{i+1}}-\max(\Delta {{\cal A}_0}_{i+1/2},0),0)/k^*_{i+1/2}\\
	  U_{i+1/2R}=(A_{i+1/2R},A_{i+1/2R}.u_{i+1})^t
        \end{array}\right.,\label{eq:rec-hydro}
\end{equation}
with $\Delta {{\cal A}_0}_{i+1/2}={{\cal A}_0}_{i+1}-{{\cal A}_0}_i=k_{i+1} \sqrt{{A_0}_{i+1}}-k_i \sqrt{{A_0}_i}$ and
 $k_{i+1/2}^*=\max(k_{i},k_{i+1})$. For consistency, the scheme \eqref{eq:scheme-hom} is modified as follows
\begin{equation}
 U_i^{n+1}=U_i^n-\dfrac{\Delta t}{\Delta x}\left(F_{i+1/2L}^n-F_{i+1/2R}^n\right),\label{eq:scheme-source}
\end{equation}
where
\begin{equation*}
 \begin{array}{l}
  F_{i+1/2L}^n=F_{i+1/2}^n+S_{i+1/2L}\\
  F_{i-1/2R}^n=F_{i-1/2}^n+S_{i-1/2R}
 \end{array},
\end{equation*}
with
\begin{equation*}
 \begin{array}{c}
  F_{i+1/2}^n={\cal F}\left(U_{i+1/2L},U_{i+1/2R},k_{i+1/2}^*\right)\\
  S_{i+1/2L}=\left(\begin{array}{c}
    0\\
    {\cal P}(A_i^n,k_i)-{\cal P}(A_{i+1/2L}^n,k_{i+1/2}^*)
  \end{array}\right)\\
 S_{i-1/2R}=\left(\begin{array}{c}
    0\\
    {\cal P}(A_i^n,k_i)-{\cal P}(A_{i-1/2R}^n,k_{i-1/2}^*)
    \end{array}\right)
 \end{array}
\end{equation*}
and ${\cal P}(A,k)=k A^{3/2}/(3\rho \sqrt{\pi})$. Thus the variation of the radius and the varying elasticity are
 treated under a well-balanced way. In system \eqref{eq:system1}, the friction term $-C_f Q/A$ is treated
 semi-implicitly. This treatment is classical in shallow water simulations \cite{Bristeau01,Liang09b} and had shown
 to be efficient in blood flow simulation as well \cite{Delestre11}. This treatment does not break the "dead man"
 equilibrium. It consists in using first \eqref{eq:scheme-source} as a prediction step without friction, {\it i.e.}:
\begin{equation*}
 U_i^*=U_i^n-\dfrac{\Delta t}{\Delta x}\left(F_{i+1/2L}^n-F_{i-1/2R}^n\right),
\end{equation*}
then we apply a semi-implicit friction correction on the predicted values ($U_i^*$):
\begin{equation*}
 A_i^*\left(\dfrac{u_i^{n+1}-u_i^*}{\Delta t}\right)=-C_f u_i^{n+1}.
\end{equation*}
Thus we get the corrected velocity $u_i^{n+1}$ and we have $A_i^{n+1}=A_i^*$.

\subsection{HLL numerical flux}\label{HLL}
As presented in \cite{Delestre11}, several numerical fluxes might be used (Rusanov, HLL, VFRoe-ncv and kinetic fluxes).
 In this work we will use the HLL flux (Harten Lax and van Leer \cite{Harten83}) because it is the best compromise
 between accuracy and CPU time consuming (see \cite[chapter 2]{Delestre10b}). It writes:
\begin{equation*}
 {\cal F}(U_L,U_R,k^*)=\left\{\begin{array}{ll}
                               F(U_L,k^*) & \text{if}\;0\leq c_1\\
			      \dfrac{c_2 F(U_L,k^*)-c_1 F(U_R,k^*)}{c_2-c_1}
+\dfrac{c_1 c_2}{c_2-c_1}(U_R-U_L) & \text{if}\;c_1<0<c_2\\
			      F(U_R,k^*) & \text{if}\;c_2\leq 0
                              \end{array}\right.,
\end{equation*}
with
\[c_{1}={\inf\limits_{U=U_L,U_R}}({\inf\limits_{j\in\{1,2\}}}\lambda_{j}(U,k^*))\;
\text{and}\;c_{2}={\sup\limits_{U=U_L,U_R}}({\sup\limits_{j\in\{1,2\}}}\lambda_{j}(U,k^*)),\label{eq:flux-hll2}\]
where \(\lambda_1(U,k^*)\) and \(\lambda_2(U,k^*)\) are the eigenvalues of the system and $k^*=\max(k_L,k_R)$.

To prevent blow up of the numerical values, we impose the following CFL (Courant, Friedrichs, Levy) condition
\begin{equation*}
 \Delta t\leq n_{CFL}\dfrac{\Delta x}{\max\limits_{i}(|u_i|+c_i)},
\end{equation*}
where $c_i=\sqrt{k_i \sqrt{A_i}/(2\rho\sqrt{\pi})}$ and $n_{CFL}=1$.

\section{Some numerical results}

\subsection{"The stented man at eternal rest"}

\begin{figure}
\begin{tabular}{cc}
  \includegraphics[angle=-90,width = 0.48 \textwidth]{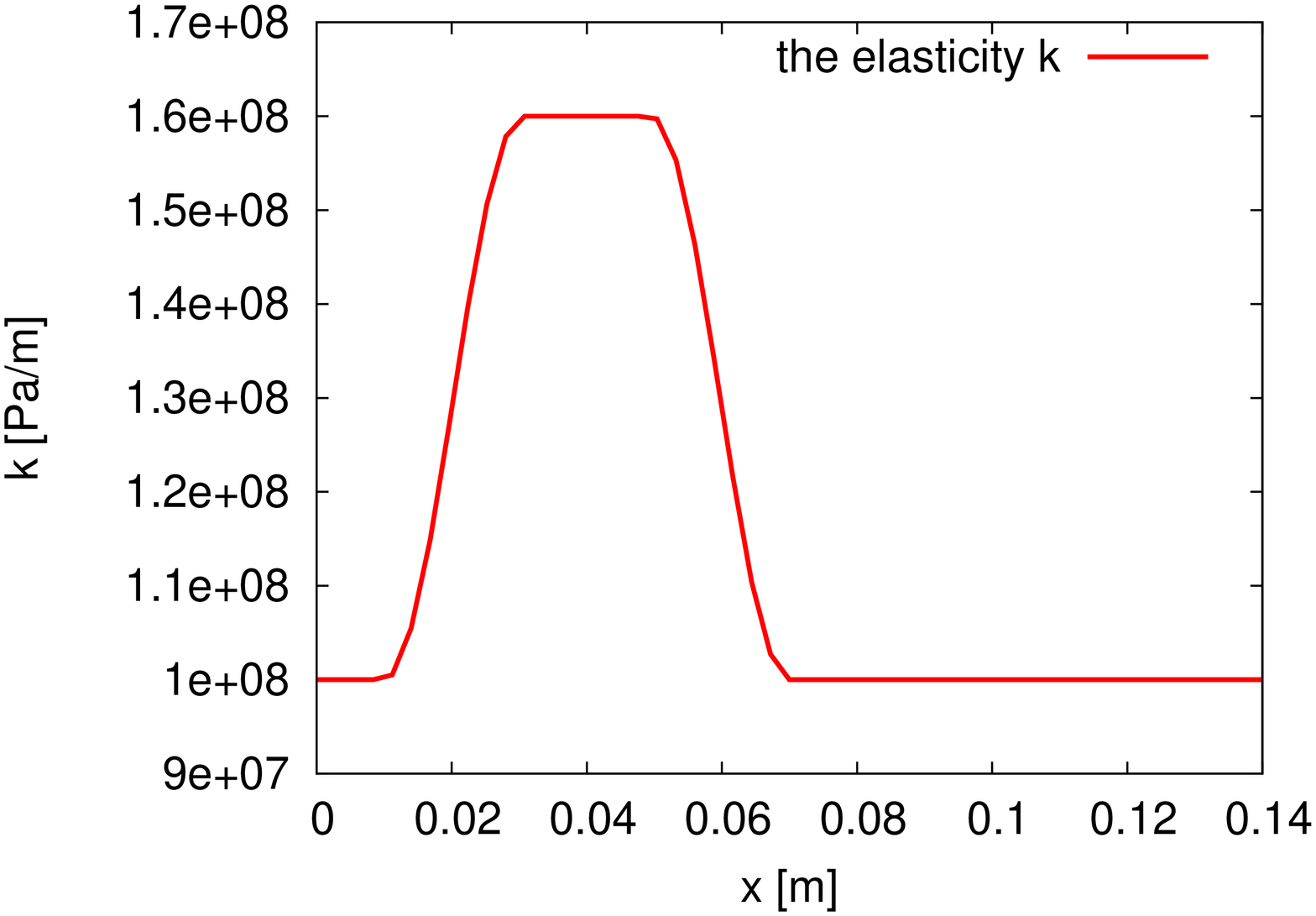} & \\[-4.2cm]
  &\includegraphics[angle=-90,width = 0.48 \textwidth]{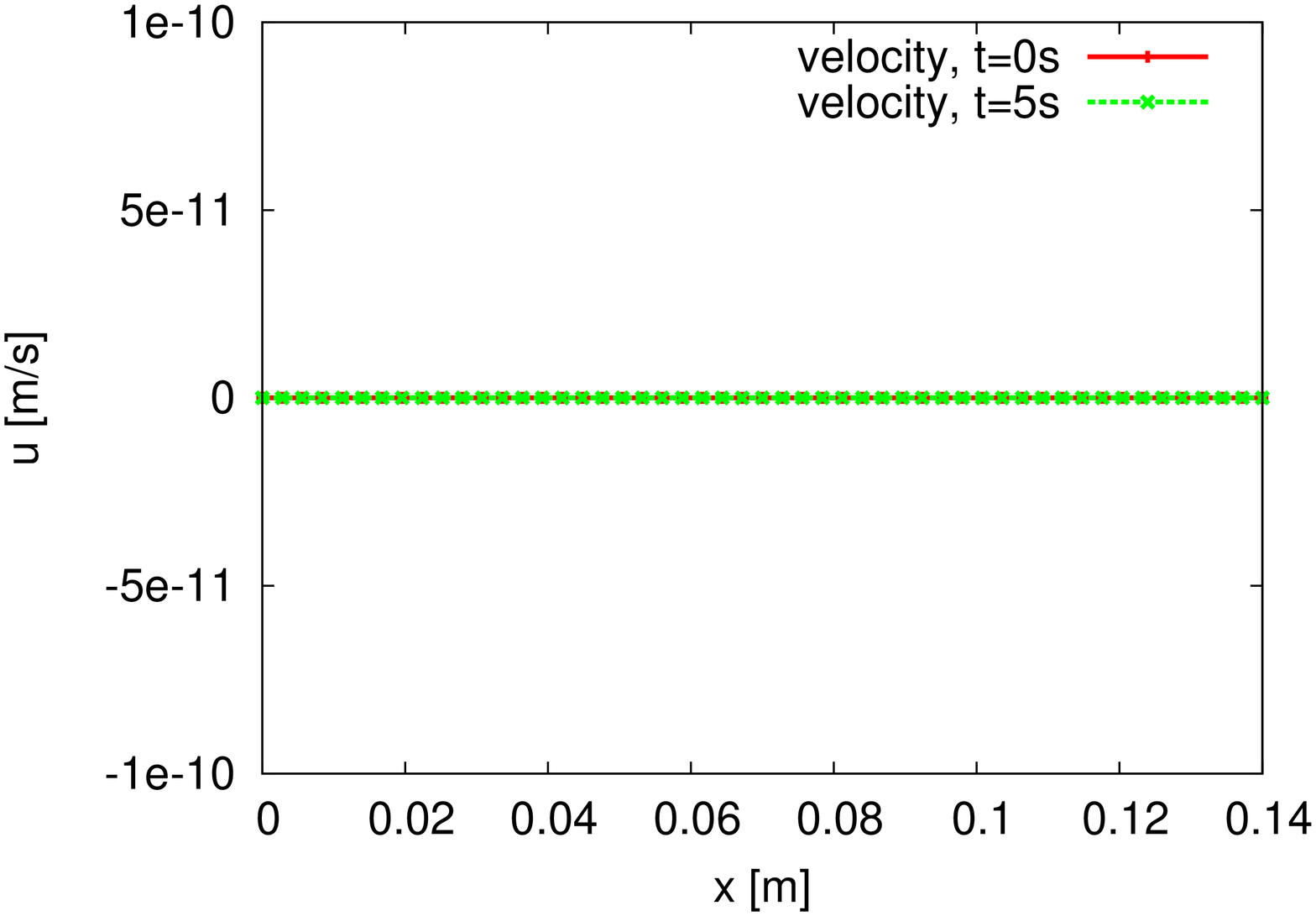}
\end{tabular}
\caption{No spurious flows (right) are generated by a change of elasticity (left).}
\label{fig:deadman}
\end{figure}

In this test, we consider a configuration with no flow and with a change of artery elasticity $k(x)$, this is
 the case for a dead man with a stented artery (see Figure 1 
 left). The section of the artery is constant
 $R_0(x)=4.0\;10^{-3}\text{m}$ and the
 velocity is $u(x,t)=0\text{m/s}$. We use the following numerical values: $J=50$ cells, $C_f=0$, $\rho=1060\text{m}^3$,
 $L=0.14\text{m}$, $T_{end}=5\text{s}$. As initial conditions, we take a fluid at rest $Q(x,0)=0\text{m}^3/\text{s}$
 and 
\begin{equation*}
 k(x)=\left\{\begin{array}{ll}
              k_0 &\text{if}\;x\in[0:x_1]\cup [x_4:L] \\
	      k_0+\dfrac{\Delta k}{2}\left(\sin\left(\dfrac{x-x_1}{x_2-x_1}\pi-\pi/2\right)+1\right) & \text{if}\;x\in]x_1:x_2[\\
	      k_0+\Delta k & \text{if}\;x\in[x_2:x_3] \\
	      k_0+\dfrac{\Delta k}{2}\left(\cos\left(\dfrac{x-x_3}{x_4-x_3}\pi\right)+1\right) & \text{if}\;x\in]x_3:x_4[
             \end{array}\right.,
\end{equation*}
with $k_0=1.0\;10^8\text{Pa}/\text{m}$, $\Delta k=6.0\;10^7\text{Pa}/\text{m}$, $x_1=1.0\;10^{-2}\text{m}$,
 $x_2=3.05\;10^{-2}\text{m}$, $x_3=4.95\;10^{-2}\text{m}$ and $x_4=7.0\;10^{-2}\text{m}$.

The steady state at rest is perfectly preserved in time, we do not notice any spurious oscillation (see Figure 1
 right).

\subsection{Wave reflection-transmission in a stented artery}

We now observe the reflexion and transmission of a pulse through a sudden change of artery elasticity (from $k_R$
 to $k_L$ with $k_L>k_R$) in an elastic tube of constant radius (see Figure 2 left). We take the following numerical
 values: $J=1500$
 cells, $C_f=0$,
 $k_L=1.6\;10^{8}\text{Pa/m}$,
 $k_R=1.\;10^{8}\text{Pa/m}$, $\Delta k=6.\;10^{7}\text{Pa/m}$,
 $\rho=1060\text{m}^3$, $R_0=4.0\;10^{-3}\text{m}$, $L=0.16\text{m}$, $T_{end}=8.0\;10^{-3}\text{s}$,
 $c_L=\sqrt{k_L R/(2\rho)}\simeq 17.37\text{m/s}$ and $c_R=\sqrt{k_R R/(2\rho)}\simeq 13.74\text{m/s}$. We take
 a decreasing elasticity on a rather small scale:
\begin{equation*}
k(x)=
\left\{\begin{array}{ll}
        k_R+\Delta k & \text{if}\; x\in[0:x_1]\\
	k_R+\dfrac{\Delta k}{2}\left[1+\cos\left(\dfrac{x-x_1}{x_2-x_1}\pi \right)\right] & \text{if}\;x\in]x_1:x_2]\\
	k_R & \text{else}
       \end{array}\right.,
\end{equation*}
with $x_1=19L/40$ and $x_2=L/2$. As initial conditions, we consider a fluid at rest $Q(x,0)=0\text{m}^3/\text{s}$ and
 the following perturbation of radius:
\begin{equation*}
 R(x,0)=\left\{\begin{array}{ll}
		  R_0(x)\left[1+\epsilon \sin\left(\dfrac{100}{20L}\pi\left(x-\dfrac{65L}{100}\right)\right)\right]
 & \text{if}\;x\in[65L/100:85L/100]\\
		  R_0(x) & \text{else}
               \end{array}\right.,
\end{equation*}
with $\epsilon=1.0\;10^{-2}$. The expression for the pressure is
\begin{equation*}
 p(x,t)=p_0(x)+k(x)(R(x,t)-R_0(x)),
\end{equation*}
where $p_0$ is the external pressure.

The numerical results perfectly match with the predictions for a linearized flow. We get the predicted amplitudes both for the
 transmitted and the reflected waves (see Figure 2 right).

\begin{figure}
\begin{tabular}{cc}
  \includegraphics[angle=-90,width = 0.48 \textwidth]{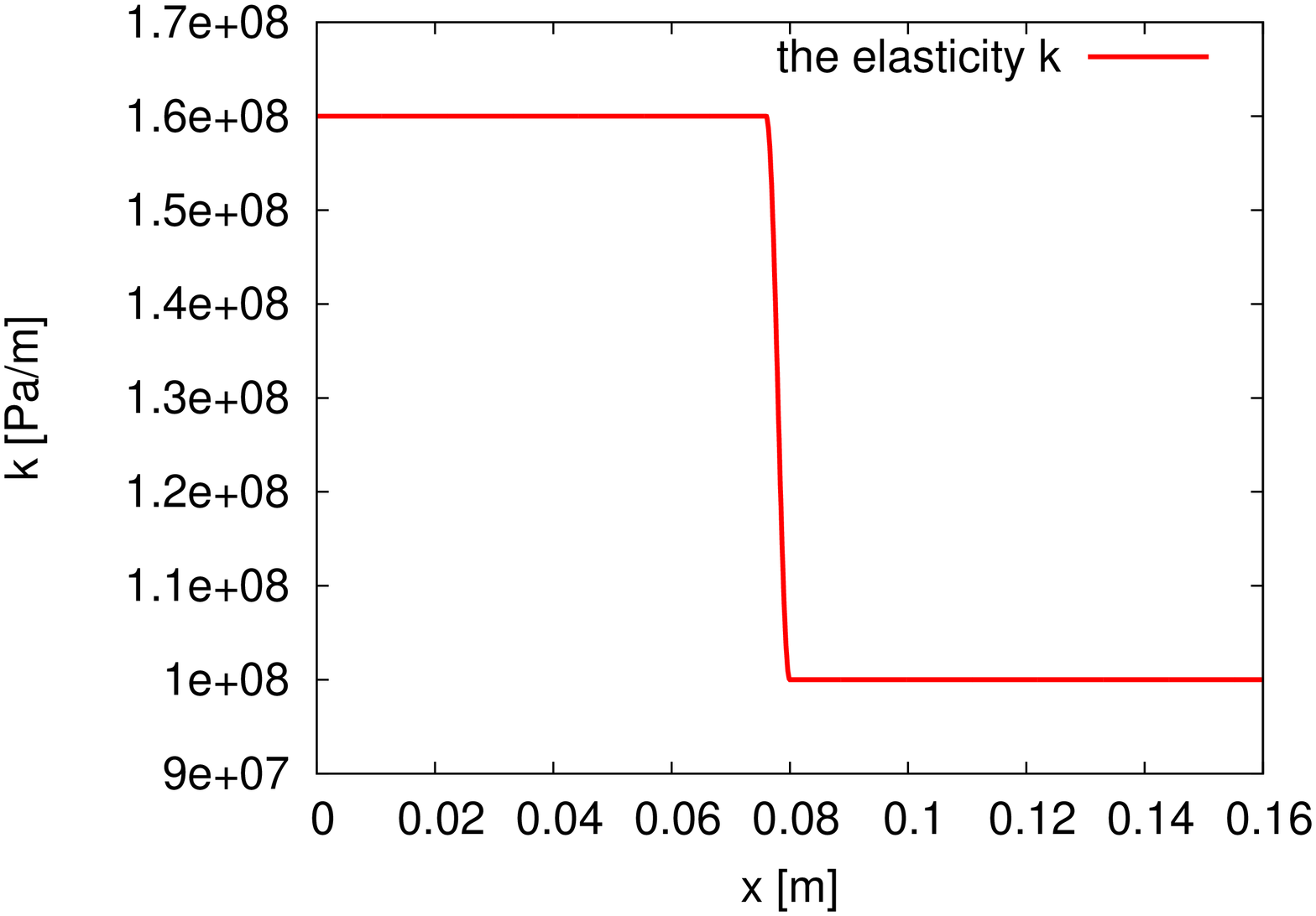} & \\[-4.2cm]
  &\includegraphics[angle=-90,width = 0.48 \textwidth]{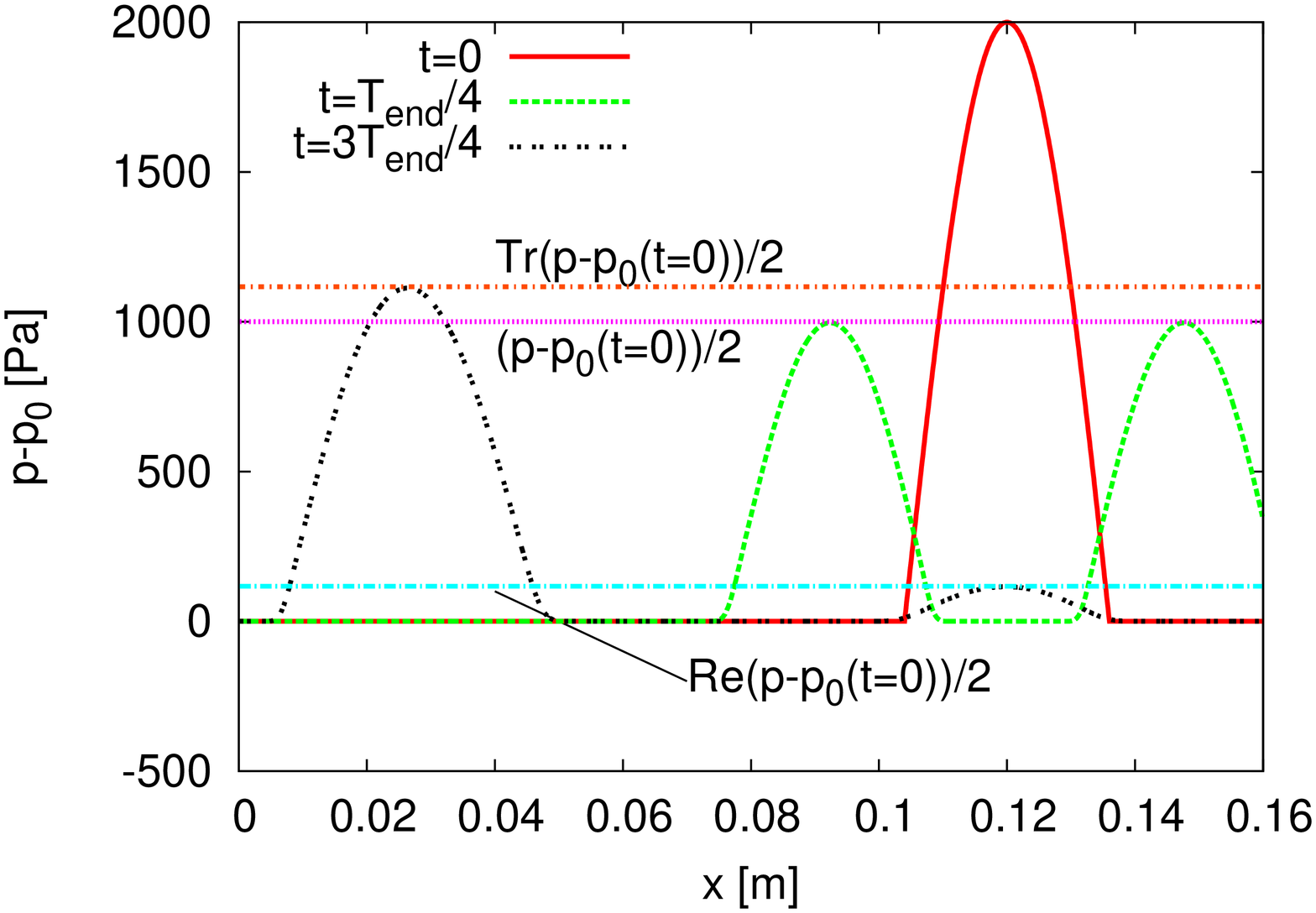}
\end{tabular}
\caption{Across a discontinuity of $k$ (left), an initial pulse evolves (right) ; it is transmitted and reflected.}
\label{fig:stented}
\end{figure}

\subsection{Wave "damping"}

In this case, the elasticity is constant in space. We consider the viscous term in the linearized momentum equation.
 A periodic signal is imposed at inflow as a perturbation of a steady state ($R_0=Cst$, $u_0=0$) with a constant
 section at rest. We take $R=R_0+\epsilon R_1$ and $u=0+\epsilon u_1$, where ($R_1$,$u_1$) is the perturbation of the
 steady state.
 Looking for progressive waves ({\it i.e.} under the form $e^{i(\omega t-Kx)}$), we take for the imposed incoming
 discharge
\begin{equation*}
 Q_b(t)=Q(t,x=0)=Q_{amp} \sin(\omega t)\;\text{m}^3/\text{s}.
\end{equation*}
Thus, we have a damping wave in the domain
 \begin{equation*}
  Q(t,x)=\left\{\begin{array}{ll}
                 0 & \text{if}\; k_r x>\omega t\\
		Q_{amp} \sin(\omega t-k_r x)e^{k_i x} &\text{if}\; k_r x \leq \omega t
                \end{array} \right.,
 \end{equation*}
with $\omega=2\pi/T_{pulse}$, $T_{pulse}$ the time length of a pulse and $K=k_r+i k_i$ the wave vector.

We use the following numerical values: $J=750$ cells, $\rho=1060\text{m}^3$, $L=3\text{m}$,
 $k=1.10^{8}\text{Pa/m}$, $R_0=4.10^{-3}\text{m}$ and $T_{end}=5\text{s}$. We consider both  $C_f=0$ and
 $C_f=0.005053$.
 As initial conditions, we take a fluid at rest
 $Q(x,0)=0\text{m}^3/\text{s}$ and as input boundary condition
\begin{equation*}
 Q_b(t)=Q_{amp}\sin(\omega t),
\end{equation*}
with $\omega=2\pi/T_{pulse}=2\pi/0.5\text{s}$ and $Q_{amp}=3.45\;10^{-7}\text{m}^3/\text{s}$. The output is an outgoing wave.

The results are closed to the analytical solution (see Figure 3). We notice a small numerical diffusion for
 $C_f=0.005053$.

\begin{figure}
\begin{tabular}{cc}
  \includegraphics[angle=-90,width = 0.48 \textwidth]{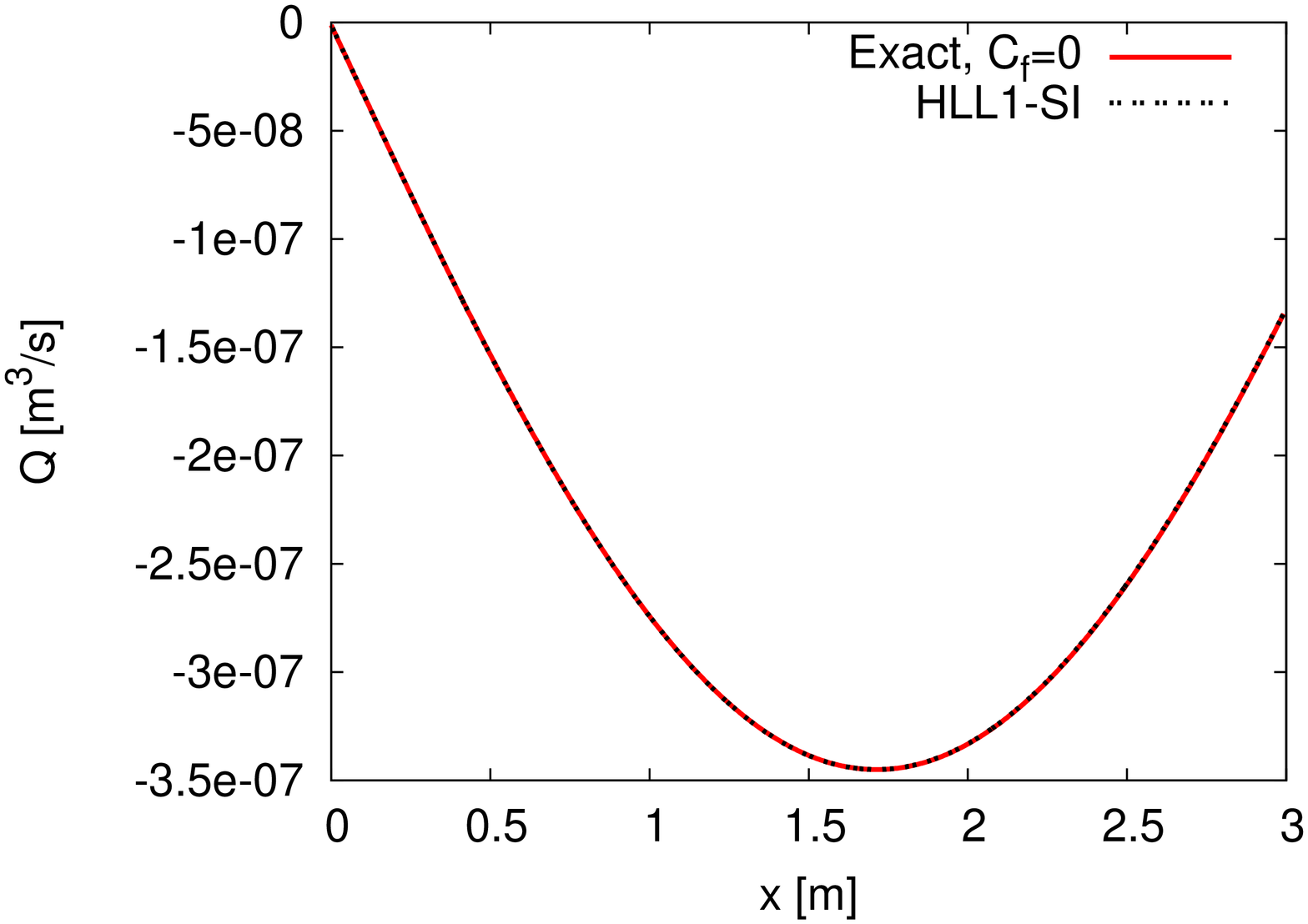} & \\[-4.2cm]
  &\includegraphics[angle=-90,width = 0.48 \textwidth]{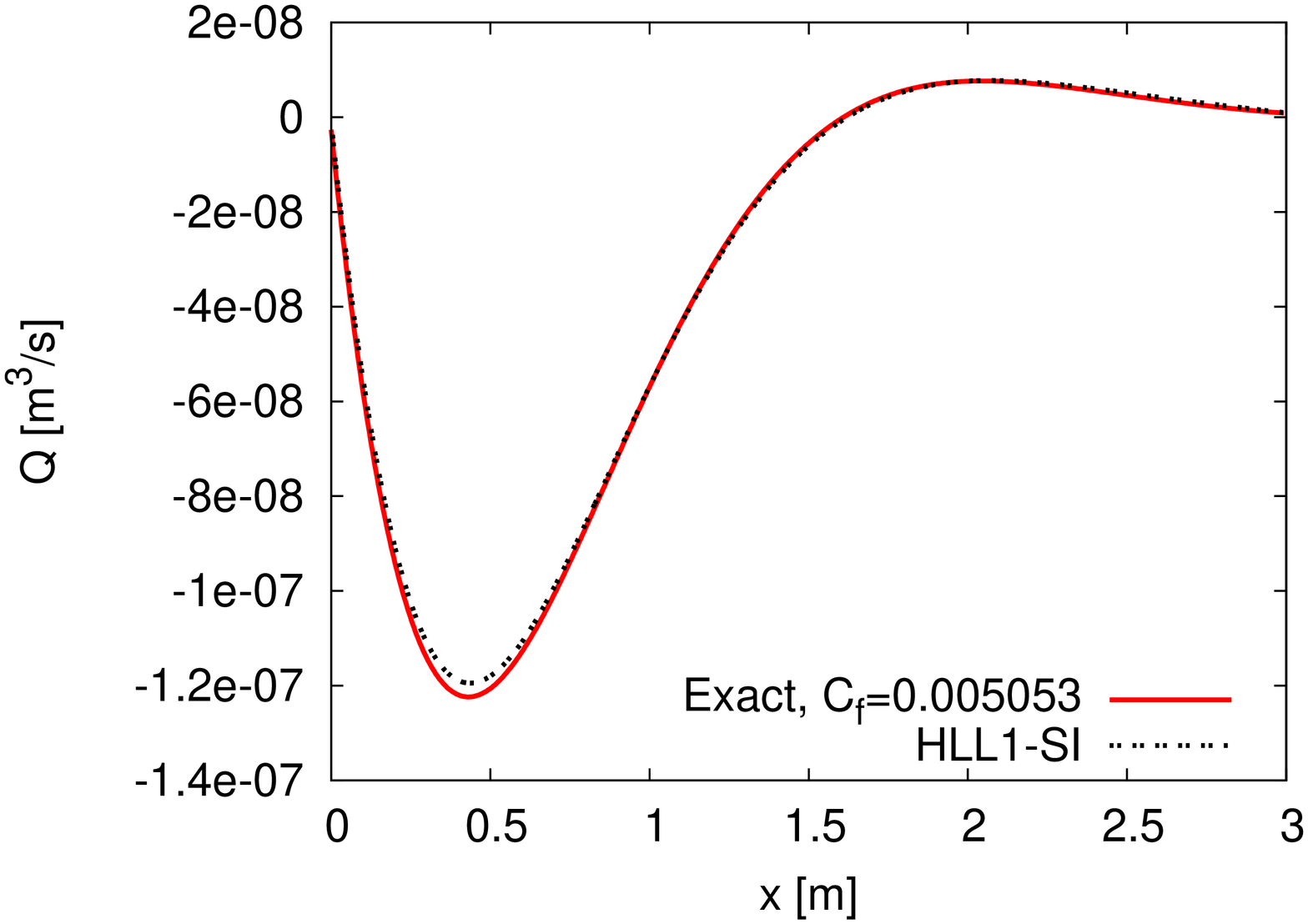}
\end{tabular}
\caption{Response to an oscillating input condition, left with no viscosity and right with viscosity.}
\label{fig::wave}
\end{figure}


\section{Conclusion}

In this work, we have considered the 1D model of flow in an artery with varying elasticity and constant section.
 We have presented a well-balanced finite volume scheme. Thus we get a mass conservative scheme. Moreover, the
 well-balanced property allows to have a good treatment of the source, {\it i.e.} we do not get numerical artifacts.
 This numerical method gave good results on numerical tests. In future works, we will have to add some extra source
 terms in order to get a more realistic model. These extra terms will require to develop a low diffusive high order
 scheme in the spirit of \cite{Delestre10}. Moreover, this will improve the accuracy of the scheme.
 And we will also have to test more complex cases such as bifurcations and networks.
\bibliographystyle{plain}

\bibliography{biblio}

\begin{thebibliography}{10}

\bibitem{Audusse04c}
E.~Audusse, F.~Bouchut, M.-O. Bristeau, R.~Klein, and B.~Perthame.
\newblock A fast and stable well-balanced scheme with hydrostatic
  reconstruction for shallow water flows.
\newblock {\em SIAM J. Sci. Comput.}, 25(6):2050--2065, 2004.

\bibitem{Bermudez94}
Alfredo Bermudez and M.~Elena Vazquez.
\newblock Upwind methods for hyperbolic conservation laws with source terms.
\newblock {\em Computers \& Fluids}, 23(8):1049 -- 1071, 1994.

\bibitem{Bouchut04}
F.~Bouchut.
\newblock {\em Nonlinear stability of finite volume methods for hyperbolic
  conservation laws, and well-balanced schemes for sources}, volume 2/2004.
\newblock Birkh{\"a}user Basel, 2004.

\bibitem{Bristeau01}
M.-O. Bristeau and Beno\^it Coussin.
\newblock Boundary conditions for the shallow water equations solved by kinetic
  schemes.
\newblock Technical Report 4282, INRIA, October 2001.

\bibitem{Delestre10b}
Olivier Delestre.
\newblock {\em Simulation du ruissellement d'eau de pluie sur des surfaces
  agricoles/ rain water overland flow on agricultural fields simulation}.
\newblock PhD thesis, Universit\'e d'Orl\'eans (in French), available from TEL:
  tel.archives-ouvertes.fr/INSMI/tel-00531377/fr, July 2010.

\bibitem{Delestre11}
Olivier Delestre and Pierre-Yves Lagr\'ee.
\newblock A "well balanced" finite volume scheme for blood flow simulation.
\newblock submitted.

\bibitem{Delestre10}
Olivier Delestre and Fabien Marche.
\newblock A numerical scheme for a viscous shallow water model with friction.
\newblock {\em J. Sci. Comput.}, DOI 10.1007/s10915-010-9393-y, 2010.

\bibitem{Greenberg96}
J.~M. Greenberg and A.-Y. Le{R}oux.
\newblock A well-balanced scheme for the numerical processing of source terms
  in hyperbolic equation.
\newblock {\em SIAM Journal on Numerical Analysis}, 33:1--16, 1996.

\bibitem{Harten83}
Amiram Harten, Peter~D. Lax, and Bram van Leer.
\newblock On upstream differencing and godunov-type schemes for hyperbolic
  conservation laws.
\newblock {\em SIAM Review}, 25(1):35--61, January 1983.

\bibitem{Kirkman03}
Robert Kirkman, Tony Moore, and Charlie Adlard.
\newblock {\em The Walking Dead}.
\newblock Image Comics, 2003.

\bibitem{Liang09b}
Qiuhua Liang and Fabien Marche.
\newblock Numerical resolution of well-balanced shallow water equations with
  complex source terms.
\newblock {\em Advances in Water Resources}, 32(6):873 -- 884, 2009.

\end{thebibliography}


\end{document}